\documentclass[12pt,a4paper,leqn, reqno]{amsart}
\usepackage{color}% insert here the call for the packages your document requires
\usepackage{latexsym}
\usepackage{amscd,amssymb,amsthm}
\usepackage{graphicx}
\def\ds{\displaystyle}
\numberwithin{equation}{section}
\theoremstyle{plain}
\newtheorem{theorem}{Theorem}[section]
\newtheorem{corollary}[theorem]{Corollary}
\newtheorem{lemma}[theorem]{Lemma}
\newtheorem{proposition}[theorem]{Proposition}

\theoremstyle{definition}
\newtheorem{definition}{Definition}[section]

\theoremstyle{remark}

\theoremstyle{remarks}

{}

\numberwithin{equation}{section}

\numberwithin{table}{section}

\numberwithin{figure}{section}

\def\ds{\displaystyle}

\newcommand{\bea}{\begin{eqnarray*}}
\newcommand{\eea}{\end{eqnarray*}}
\newcommand{\bean}{\begin{eqnarray}}
\newcommand{\eean}{\end{eqnarray}}

\newcommand{\al}{\alpha}

\newcommand{\la}{\lambda}

\begin{document}
%Paley type inequality for  Heisenberg group $H^p(\mathbb{H}^{n})$
% ----------------------------------------------------------------
\title[Paley type inequality for the Fourier transform on the Heisenberg group]
{\bf Paley type inequality  of the Fourier transform on the Heisenberg group}
\author{\bf Atef Rahmouni}
\address{Department of mathematics, King Saudi University, College of Sciences\\
P. O Box 2455 Riyadh 11451, Saudi Arabia.}
\email{Atef.Rahmouni@fsb.rnu.tn}

% ----------------------------------------------------------------
\begin{abstract}
A paley type inequality  for the Fourier transform on  $H^p(\mathbb{H}^{n}),$ the Hardy space on the Heisenberg group, is obtained for $0<p\leq 1.$
\end{abstract}

% ----------------------------------------------------------------
\maketitle
% ----------------------------------------------------------------
\section{ Introduction}

The study of Hardy spaces has been  originated during the 1910's
in the setting of Fourier series and complex analysis in one
variable. In 1972, Fefferman and Stein \cite{FS} introduced Hardy
spaces $H^p$ by mean of maximal function
$$f^\ast(x) = \sup_{r>0} |f\ast\phi_r(x)|$$
where  $\phi $ belongs to $\mathcal{S}$, the Schwartz space of
rapidly decreasing smooth functions  satisfying $\int \phi(x) dx =
1$. The delation   $\phi_r$ is given by  $\phi_r(x) =
r^{-n}\phi(x/r).$
 We say that a tempered distributions $f\in\mathcal{ S}'$ is in $H^p$ if $f^\ast$
is in $L^p$. \\ Using the maximal function above, Coifman
\cite{CO} showed that any $f$ in $H^p$ can be represented as a
linear combination of atoms, that is
$$ f = \sum_{k=1}^{\infty}\beta_{k}a_k,\quad \beta_{k}\in \mathbb{C},
$$ where the $a_k$ are atoms and the sum converges in $H^p$. Moreover,
$$\|f\|_{H^p} \thickapprox \inf \Big\{ \sum_{k=1}^\infty |\beta_{k}|^p :
\sum_{k=1}^{\infty}\beta_{k}a_{k} \mbox{ is a decomposition of $f$
into\, }  atoms\Big\}.$$ It has been shown that  the study of some
analytic problems on $H^p(\mathbb{R}^n)$ is summed up to
investigate some properties of these atoms, and therefore the
problems become quite simple. In 1980, Taibleson and Weiss
\cite{TW} gave the definition of molecules belonging to $H^p,$ and
showed that every molecule is in $H^p$ with continuous embedding
map. By the atomic decomposition and the molecule
characterization, the proof of $H^p$ boundedness of the operators
on Hardy space becomes easier. The theory of $H^p$ have been
extensively studied in \cite{GR} and \cite{S3}.

In the setting of the euclidian case, Hardy's inequality for
Fourier transform asserts that for all $f\in H^p(\mathbb{R}^n)$
$0<p\leq 1.$
\begin{equation} \label{A}
\int_{\mathbb{R}^n}{|\widehat{f}(\xi)|^p\over |\xi|^{n(2-p)}}
d\xi\leq \|f\|^p_{H^p{(\mathbb{R}^n})},\qquad 0<p\leq1
\end{equation}
where $H^p(\mathbb{R}^n)$ indicates the real Hardy space.  Hardy's
type inequality for Fourier transform has been extensively studied
in \cite{S2}. Kanjin \cite{kan} proved Hardy's inequalities for
Hermit and Laguerre expansions for functions in $H^1$  and for
Hankel transform \cite{kan2}. In connection with properties of regularity of the spherical means on $\mathbb{C}^n$, Thangavelu \cite{thang} proved
a Hardy's inequality for special Hermit functions. These standard
inequalities for higher dimensional has been  studied in \cite{RT}. Recently, an extension has been given
by \cite{AM1}, the latter establish a Hardy's type inequality
associated with the Hankel transform for over critical exponent
$\sigma >\sigma_0= 2 - p.$ We point out here that the result
obtained for Hardy's inequality for the Hankel transform improves
the work of Kanjin \cite{kan2} in which he proved the result for $
\sigma_0=2 - p.$ Although, in \cite{AR,AR1,AR3}  extended this form
of this inequality to Laguerre hypergroup and its dual.

In this paper we are interested in the Heisenberg group $\mathbb{H}^n$ is the lie group with underlying manifold $\mathbb{H}^n=\mathbb{C}^n \times \mathbb{R}$ and multiplication $(z,t).(z',t')=(z+z',t+t'+2 Im(z.\overline{z'}),$ where $\linebreak z=(z_1,z_2,...,z_n)\in \mathbb{C}^n.$ If we identify $\mathbb{C}^n \times \mathbb{R}$ with $\mathbb{R}^{2n+1}$ by $z_j=x_j+ix_{j+n},~j=1,...,n,$
then the group law can be rewritten as $$(x_1,x_2,...,x_{2n},t).(y_1,y_2,...,y_{2n},t')=(x_1+y_2,...,
x_n+y_n,t+t'-2\sum_{j=1}^{n}(x_jy_{j+n}-y_jx_{j+n})).$$
The reverse element of $u=(z, t)$ is $u^{-1}=(-z,-t)$ and we write the identity of $\mathbb{H}^n$ as $0=(0,0).$

Set $X_{j}, X_{j+n}$ and $T$ is a basis for the left invariant vector fields on $\mathbb{H}^n.$
The corresponding complex vector fields are
$$Z_j=\frac{1}{2}(X_{j}-iX_{j+n})=\frac{\partial}{\partial z_{j}}+i\overline{z}_{j}\frac{\partial}{\partial t},~
\overline{Z}_j=\frac{1}{2}(X_{j}+iX_{j+n})=\frac{\partial}{\partial \overline{z}_{j}}-iz_{j}\frac{\partial}{\partial t},~j=1,...,n.$$
The Heisenberg group is a connected, simply connected nilpotent Lie group. We define one-parameter dilations on $\mathbb{H}^n,$  for $R>0,$ by
$\rho_{R}(z,t)=(R z,R^2 t).$ These dilations are group automorphisms and the Jacobian determinant is $R^{Q},$ where $Q = 2n + 2$ is the \textit{homogeneous dimension} of $\mathbb{H}^n.$ We will denoted by
$f_{\rho}(z,t) =\rho^{-Q}f((z,t)_{\rho})$ the dilated of the
function $f$ defined on $\mathbb{H}^{n}.$

A homogeneous norm on $\mathbb{H}^{n}$ is given by
$$|(z, t)|_{\mathbb{H}^{n}}=(|z|^{4}+4t^{2})^{1/4},$$
With this norm, we define the Heisenberg ball centered at $u=(z, t)$ of radius $r,$ i.e., the set
$$B(u,r)=\{v\in\mathbb{H}^{n}:~|uv^{-1}|_{\mathbb{H}^{n}}<R\},$$
and we denote by $B_R = B(0, R) = \{v\in\mathbb{H}^{n}:~|v|_{\mathbb{H}^{n}}<R\} $ the open ball centered at 0, the identity element of $\mathbb{H}^{n},$ with radius $R.$ The volume of the ball $B(u, R)$ is $C_QR^Q,$ where $C_Q$ is the volume of the unit ball $B_1.$

The Haar measure $dV$ on $\mathbb{H}^{n}$ coincides with the Lebesgue measure on $\mathbb{C}^n \times \mathbb{R}$ which is denoted by $dzd\overline{z}dt.$

Let $J=(j^{1},j^{2},j^{0})\in
\mathbb{Z}^n_{+} \times \mathbb{Z}^n_{+}\times \mathbb{Z}_{+},$ where $\mathbb{Z}_{+}$ the set of all nonnegative integers, we set $h(J)=|j^{1}|+|j^{2}|+2j^{0},$
where, if $j^{1}=(j^{1}_{1},...,j^{n}_{n}),$ then
$|j^1|=\sum_{k=1}^{n}j^{1}_{k}.$ If $P(z,t)=\sum_{J}a_{J}(z,t)^{J} $ is a polynomial
where $(z, t)^J = z^{j^{1}}\overline{z}^{j^{2}}t^{j^{0}},$ then we call $\max\{h(J) : a_J \neq 0\}$ the homogeneous
degree of $P(z, t).$ The set of all polynomials whose homogeneous degree $\leq s$ is denoted by $\mathcal{P}_s.$ Schwartz space on $\mathbb{H}^{n}$ write as $\mathcal{S}(\mathbb{H}^{n}).$

Fix $\lambda > 0,$ let $\mathcal{H}_{\lambda}$ be the Bargmann's space :
$$\mathcal{H}_{\lambda}=\Big\{F ~\mbox{holomorphic on}~\mathbb{C}^{n}: \|F\|^{2}=\big(\frac{2\lambda}{\pi}\big)^{n}\int_{\mathbb{C}^{n}}
|F(\zeta)|^{2}e^{-2\lambda|\zeta|^{2}}d\zeta<\infty \Big\}.$$
Then, $\mathcal{H}_{\lambda}$ is a Hilbert space and the monomials
$$F_{\alpha,\lambda}(\zeta)=\sqrt{\frac{(2\lambda)^{|\alpha|}
}{\alpha!}}\zeta^{\alpha},~~~\alpha=(\alpha_{1},\alpha_{2},...,
\alpha_{n})\in
\mathbb{Z}^{n}_{+}$$
form an orthonormal basis for $\mathcal{H}_{\lambda},$ where $\alpha! = \alpha_{1}!\alpha_{2}!...\alpha_{n}!,~ |\alpha| = (\alpha_{1},\alpha_{2},...,\alpha_{n})$ and $\linebreak\zeta^{\alpha}= \zeta^{\alpha_{1}}_{1}\zeta^{\alpha_{2}}_{2}...
\zeta^{\alpha_{n}}_{n}.$ Suppose $W_{k,\lambda}$ and
$W^{+}_{k,\lambda}$ are the closed operators on
$\mathcal{H}_{\lambda}$ such that
\begin{eqnarray*}
W_{k,\lambda}F_{\alpha,\lambda}&=&(2(\alpha_{k}+1)\lambda)^{1/2}
F_{\alpha+e_{k},\lambda},\\
W_{k,\lambda}^{+}F_{\alpha,\lambda}&=&(2\alpha_{k}\lambda)^{1/2}
F_{\alpha-e_{k},\lambda},\qquad \mbox{for}~~\lambda>0,
\end{eqnarray*}
and
\begin{eqnarray*}
W_{k,\lambda}&=&W_{k,-\lambda}^{+},\\
W_{k,\lambda}^{+}&=&W_{k,-\lambda},\qquad \mbox{for}~~\lambda<0,
\end{eqnarray*}
where $e_k = (0,...,1,...,0)\in\mathbb{Z}^n$ with the 1 in the $k$-th position. Then
$$\prod_{\lambda}(z,t)=exp^{i\lambda t}exp^{(-z.W_{\lambda}
+\overline{z}.W_{\lambda}^{+})}$$
is an irreducible unitary representation of $\mathbb{H}^n$ on $\mathcal{H}_{\lambda},$ where $z.W_{\lambda}=\sum_{k=1}^{n}
z_{k}.W_{k,\lambda}.$

The group Fourier transform of $f\in L^1(\mathbb{H}^n)\cap L^2(\mathbb{H}^n)$ is an operator-valued function
defined by
\begin{equation}
\mathcal{F}(f)(\lambda)=\int_{\mathbb{H}^n}f(z,t)
\prod_{\lambda}(z,t)dV.
\end{equation}
Obviously, $\|\mathcal{F}(f)(\lambda)\|\leq \|f\|_{L^{1}}.$
Here, $\|-\|$ denotes the operator norm. Similar as in
$\mathbb{R}^n,$ for $f\in L^1(\mathbb{H}^n)\cap L^2(\mathbb{H}^n),$ we have the following Plancherel and inversion formulas :
\begin{equation}
\|f\|_{2}^{2}=\frac{2^{n-1}}{\pi^{n+1}}\int_{\mathbb{R}}
\|\mathcal{F}(f)(\lambda)\|^{2}_{HS}|\lambda|^{n}d\lambda,
\qquad f\in L^1(\mathbb{H}^n)\cap L^2(\mathbb{H}^n),
\end{equation}
\begin{equation}
\int_{\mathbb{R}}tr\Big(\prod^{*}_{\lambda}(z,t)
\mathcal{F}(f)(\lambda)\Big)|\lambda|^{n}d\lambda=
\frac{(2\pi)^{n+1}}{4^n} f(u)
\end{equation}
where $tr$ is the canonical semifinite trace
and $\|-\|_{HS}$ denotes the Hilbert-Schmidt norm.\\
For $(\lambda,m, \alpha) \in\mathbb{R}^{*} \times\mathbb{Z}^{n}\times\mathbb{Z}^{n}_{+},$
where $\mathbb{R}^{*}=\mathbb{R}\backslash\{0\},$ we use the notations
\begin{center}
\begin{tabular}{cll}
 & $m_{i}^{+}=\max\{m_{i},0\},\qquad$ & $m_{i}^{-}=-\min\{m_{i},0\},$\\
&$m^{+}=(m_{1}^{+},m_{2}^{+},...,m_{n}^{+}) \qquad$ &$m^{-}=(m_{1}^{-},m_{2}^{-},...,m_{n}^{-}).$
\end{tabular}
\end{center}
The partial isometry operator $W^{m}_{\alpha}(\lambda)$
on $\mathcal{H}_{|\lambda|}$  by
\begin{center}
\begin{tabular}{cll}
 & $W_{k,\alpha}(\lambda)F_{\beta,\lambda}=(-1)^{|m^{+}|}
\delta_{\alpha+m^{+},\beta}F_{\alpha+m^{-},\lambda},$
&$\qquad\mbox{for}~~\lambda>0;$\\
& $W_{\alpha}^{m}(\lambda)=[W_{\alpha}^{m}(-\lambda)]^{*},$
&$\qquad    \mbox{for}~~\lambda<0.$
\end{tabular}
\end{center}
Thus $\{W^{m}_{\alpha}(\lambda): m\in\mathbb{Z}^n,\alpha\in\mathbb{Z}^n\}$ is an orthonormal basis for the Hilbert-Schmidt operators on $\mathcal{H}_{|\lambda|}.$ Given a function $f\in L^2(\mathbb{H}^n)$ such that
$$f(z,t)=\sum_{m,\alpha}f_{m}(r_{1},...,r_{n},t)
e^{i(m_{1}\theta_{1}+...+m_{n}\theta_{n})},\qquad\mbox{where}
\qquad z_{j}=r_{j}e^{i\theta_{j}},$$
then
$$\mathcal{F}(f)(\lambda)=\sum_{m,\alpha}R_{f}(\lambda,m,\alpha)
W^{m}_{\alpha}(\lambda),$$
where
$$R_{f}(\lambda,m,\alpha)=\int_{\mathbb{H}^{n}}
f_{m}(r_{1},...,r_{n},t)e^{i\lambda t} \ell_{\alpha_{1}}^{|m_{1}|}(2|\lambda|r_{1}^{2})...
\ell_{\alpha_{n}}^{|m_{n}|}(2|\lambda|r_{n}^{2})dV,$$
and $\ell_{\alpha}^{|m|}$ is the Larguerre function of type $|m|$ and degree $|\alpha|.$

Let $P$ be a polynomial in $z_j, \overline{z}_j, t$ on $\mathbb{H}^n,$ and we define the difference-differential
operator $\Delta_P$ acting on the Fourier transform of
$f\in L^1 \cap L^2(\mathbb{H}^n)$ by
$$\Delta_P\Big(\sum_{m,\alpha}R_{f}(\lambda,m,\alpha)
W_{\alpha}^{m}(\lambda)\Big)=\sum_{m,\alpha}R_{Pf}
(\lambda,m,\alpha)W_{\alpha}^{m}(\lambda),$$
namely, $\Delta_P\mathcal{F}(f)(\lambda)=\widehat{P(.)f(.)}(\lambda).$
In \cite{Liu} and \cite{Lin}, the authors gave the explicit expressions
for $\Delta_t,\Delta_{z_{j}}$ and $\Delta_{\overline{z}_j}.$ For convenience, we shall write $\Delta^{J}_{(z,t)} = \Delta^J.$

The paper is organized as follows. In the Second section we give an appropriate definition of atoms and investigate the atoms characterization of Hardy spaces
$H^p(\mathbb{H}^{n})$ for $0 <p\leq 1.$  In the last section  we
state  and prove our main  result:
\begin{theorem}\label{t3.3.11}
Let $0<p \leq 1,$ and  $ s\geq J=[Q({1/p}-1)], $ the greatest integer
not exceeding $\linebreak Q({1/p}-1).$ Then for any $f\in H^{p}(\mathbb{H}^{n})$ the Fourier transform of $f$ satisfies
the following Hardy's type inequality
\begin{equation}\label{3.3.13}
\int_{\mathbb{R}}\frac{\|\mathcal{F}(f)(\lambda)\|_{HS}^{p}}
{\big((2|\alpha|+n)|\lambda|\big)^{\sigma}}
|\lambda|^{n}d\lambda\leq C(p,n)\|f\|^{p}_{H^{p}(\mathbb{H}^{n})},
\end{equation}
provided that
\begin{equation}\label{3}
\frac{Q}{2}(2-p)\leq \sigma<\frac{Q}{2}+p(\frac{J+1}{2})
\end{equation}
where $C(p, n)$ depend only on $p$ and $n.$
\end{theorem}
Finally, we mention that $C$ will be always used to denote a
suitable positive constant that is not necessarily the same in
each occurrence.

\section{Atomic decomposition for $H^{p}(\mathbb{H}^{n})$}

Now we state the definition of atomic Hardy spaces in the setting
of the Heisenberg group  $H^p(\mathbb{H}^{n})$,\ $0<p\leq1$. To
this end, we introduce the following kind of atoms, which is
closely related to the Haar measure $dV.$
\begin{definition}
Let $0 < p \leq 1 \leq q \leq \infty,  p \neq q,s\in\mathbb{Z}$ and $ s \geq J=[Q(1/p- 1)].$ (Such an ordered triple $(p, q, s)$ is called admissible).
A $(p, q, s)$-atom centered at $x_{0} \in \mathbb{H}^{n}$ is a function $a \in L^{q}(\mathbb{H}^{n}),$ supported on a ball
$B(x_{0},R)\subset \mathbb{H}^{n}$ with centre $x_0=(z_0,t_0)$ and
satisfying the following
\begin{itemize}
\item[(i)] $\|a\|_{L^q(\mathbb{H}^{n})}\leq |B(0,r)|^{\frac{1}{q}-{1\over p}} $,
a.e,
\item[(ii)] $\ds \int_{\mathbb{H}^{n}} a(x)P(x)dV(x)= 0 ,$ for every $P\in \mathcal{P}_{s}.$
\end{itemize}
\end{definition}
Here, $(i)$ means that the size condition of atoms, and $(ii)$ is
called the cancelation moment condition.

A characterization of $H^p(\mathbb{H}^{n})$ is included in the following statements.
\begin{proposition}
Let $0 < p \leq1.$ If $\{a_{k} \}_{k=0}^{\infty}$ is a sequence of $p$-atoms, and $\{\lambda_k \}_{k=0}^{\infty}$ is a sequence of complex numbers with
$$\Big(\sum_{k=0}^{\infty} |\lambda_{k}|^{p}\Big)^{1/p}<\infty,$$
then $\sum_{k=0}^{\infty} \lambda_{k}a_{k} $ converges in $H^p(\mathbb{H}^{n})$ and
$$\Big\| \sum_{k} \lambda_{k}a_{k} \Big\|_{H^p(\mathbb{H}^{n})}\leq C(p,n) \Big(\sum_{k} |\lambda_{k}|^{p}\Big)^{1/p}  .$$

Conversely, if $f\in H^p(\mathbb{H}^{n})$ there exists a sequence $\{a_k \}_{k=0}^{\infty}$ of $p$-atoms, and a sequence $\{\lambda_k \}_{k=0}^{\infty}$ of complex numbers such that
$$ f=\sum_{k} \lambda_{k}a_{k}\qquad\mbox{and}\qquad
\Big(\sum_{k} |\lambda_{k}|^{p}\Big)^{1/p}\leq  C(p,n) \|f\|_{H^p(\mathbb{H}^{n})},$$
where $C(p,n)$ depends on $p$ and $n.$
\end{proposition}

\section{Proof of the main  result}

Now we are in a position to give the proof of the  main  result.
First we stat  the following proposition which has its own interest.

\begin{proposition}\label{p3.3.10}
For all $(z,t)\in \mathbb{H}^{n}$  the function
$\prod_{\lambda}(z,t)$ satisfies

\begin{equation}
\prod_{\lambda}(z,t)=\sum_{2k+\ell\leq J}
\omega_{k,\ell}(\lambda, n)~ z^{k}t^{\ell}+R_{\theta}(z,t), \,\,
 0<\theta<1,
\end{equation}
where
\begin{equation}
R_{\theta}(z,t)=\sum_{2k+\ell=J+1}
\frac{(i \lambda t)^{k}} {k!}.\frac{ (z.W_{\lambda}-\overline{z}.W^{+}_{\lambda})^{\ell}}{\ell!}.
%|R_{\theta}(z,t)|\leq C~\sum_{2k+\ell=
%J+1}z^{\ell}|t|^{k}[\mathcal{N}(\lambda,n)]^{\frac{k}{2}+\ell}, \quad
%\mbox{if } \, i \quad  \mbox{is even}.
\end{equation}
Here  $ \omega_{k,\ell}(\lambda, n)$ are functions expressed by mean of $\lambda, n.$

Set $\mathcal{H}_{|\lambda|}^{N}$ be the subspace of $\mathcal{H}_{|\lambda|}$ spanned by $\{W^{0}_{\alpha}(\lambda) :|\alpha|\leq N\}.$
Remark that (see \cite{Liu,Lin1}) $z.W_{\lambda}-\overline{z}.W^{+}_{\lambda}$ is
bounded from $\mathcal{H}_{|\lambda|}^{N}$  to $\mathcal{H}_{|\lambda|}^{N+1}$
and whose bound $< ((2|\alpha| +n)|\lambda|)^{1/2}|z|.$
Then
$$R_{\theta}(z,t)\leq C \sum_{2k+\ell=J+1}\omega_{k,\ell}~
\big((2|\alpha|+n)|\lambda|\big)^{k+\frac{\ell}{2}}~z^{\ell}~t^{k}.$$
\end{proposition}
{\bf Proof of Theorem \ref{t3.3.11}.} Let
$f=\sum_{k=0}^{\infty}\beta_{k}a_{k}\in H^{p}(\mathbb{H}^{n}),$ being
element of $H^{p}(\mathbb{H}^{n})$ where $a_{k}$ are atoms. Since $0 < p \leq 1 $ it follows
$$\int_{\mathbb{R}}\frac{\|\mathcal{F}(f)(\lambda)\|_{HS}^{p}}
{\big((2|\alpha|+n)|\lambda|\big)^{\sigma}}
|\lambda|^{n}d\lambda\leq C
\sum_{k=0}^{\infty}|\beta_{k}|^{p}\int_{\mathbb{R}}
\frac{\|\mathcal{F}(a_{k})(\lambda)\|_{HS}^{p}}
{\big((2|\alpha|+n)|\lambda|\big)^{\sigma}}
|\lambda|^{n}d\lambda.$$ In order to prove  Theorem \ref{t3.3.11}, it is enough to prove,
\begin{equation}\label{3.3.15}
\int_{\mathbb{R}}\frac{\|\mathcal{F}(a_{k})(\lambda)\|_{HS}^{p}}
{\big((2|\alpha|+n)|\lambda|\big)^{\sigma}}
|\lambda|^{n}d\lambda\leq C.
\end{equation}
This follows as $f=\sum_{k=0}^{\infty}\beta_{k}a_{k}$ implies
$\mathcal{F}(a_{k})(\lambda)^{p}\leq\Big|\sum_{k}\beta_{k}
\mathcal{F}(a_{k})(\lambda)\Big|^{p} \leq \sum_{k=0}^{\infty}|\beta_{k}|^{p}|\mathcal{F}(a_{k})(\lambda)|^{p}$
and hence
\begin{eqnarray*}
\int_{\mathbb{R}}\frac{\|\mathcal{F}(f)(\lambda)\|_{HS}^{p}}
{\big((2|\alpha|+n)|\lambda|\big)^{\sigma}}
|\lambda|^{n}d\lambda&\leq & C
\sum_{k=0}^{\infty}|\beta_{k}|^{p}
\int_{\mathbb{R}}\frac{\|\mathcal{F}(a_{k})(\lambda)\|_{HS}^{p}}
{\big((2|\alpha|+n)|\lambda|\big)^{\sigma}}
|\lambda|^{n}d\lambda\\
&\leq& C \Big\{\sum_{k=0}^{\infty}|\beta_{k}|^{p}\Big\}^{1/p}\\
&\leq & C\|f\|_{H^{p}(\mathbb{H}^{n})}.
\end{eqnarray*}
Let us  now take $\gamma$ an arbitrary nonnegative real number,
and decomposing the left hand side of (\ref{3.3.15}) as {\small
\begin{eqnarray*}
\int_{\mathbb{R}}\frac{\|\mathcal{F}(a_{k})(\lambda)\|_{HS}^{p}}
{\big((2|\alpha|+n)|\lambda|\big)^{\sigma}}
|\lambda|^{n}d\lambda
\!\!\!\!\!\!\!&\!\!\!\!\!=\!\!\!\!\!\!\!\!\!&\!\!\!
\int_{0<|\lambda|\leq\gamma}\frac{\|\mathcal{F}(a_{k})(\lambda)
\|_{HS}^{p}}
{\big((2|\alpha|+n)|\lambda|\big)^{\sigma}}|\lambda|^{n}d\lambda\\\nonumber\\\\
&+&
\int_{|\lambda|>\gamma}\frac{\|\mathcal{F}(a_{k})(\lambda)\|_{HS}^{p}}
{\big((2|\alpha|+n)|\lambda|\big)^{\sigma}}
|\lambda|^{n}d\lambda\nonumber\\\nonumber\\\\
&\,\,\,\,:=& S_{1}+S_{2}.
\end{eqnarray*}
To estimate $S_{1}$ we may use Proposition \ref{p3.3.10}, and
cancelation property of atoms. Hence, by the cancelation property of atom,
$$\mathcal{F}(a_{k})(\lambda)=\int_{\mathbb{H}^{n}}\Big[\sum_{2k+\ell\leq J}
\omega_{k,\ell}(\lambda, n)~ z^{k}t^{\ell}+R_{\theta}(z,t)\Big]
a (z,t)~dV(z,t).$$ Now with the help of properties $(i),
(ii)$ for $a(p,\infty,s)$-atoms of $H^{p}(\mathbb{H}^{n})$ together
with  Proposition \ref{p3.3.10}, we get
\begin{eqnarray*}
\mathcal{F}(a_{k})(\lambda)&\leq &
C \sum_{2k+\ell=J+1}\omega_{k,\ell}~
\big((2|\alpha|+n)|\lambda|\big)^{k+\frac{\ell}{2}}
\int_{B(o,R)}
~z^{\ell}~t^{k}
|B(0,R)|^{-\frac{1}{p}}~dV(z,t)\nonumber\\
&\leq& C ~ \sum_{2k+\ell= J+1}\omega_{k,\ell}
~R^{Q(1-\frac{1}{p})+2k+\frac{\ell}{2}} ~\big((2|\alpha|+n)|\lambda|\big)^{k+\frac{\ell}{2}}.
\end{eqnarray*}
Integrating with respect to the measure
$d\gamma_{n}(\lambda)=|\lambda|^{n}d\lambda$ over the domain $0\leq
|\lambda|\leq \gamma,$ we obtain {\small
\begin{eqnarray*}
S_{1}&= &
\int_{0<|\lambda|\leq\gamma}\frac{\|\mathcal{F}(a_{k})(\lambda)\|^{p}}
{\big((2|\alpha|+n)|\lambda|\big)^{\sigma}}|\lambda|^{n}
d\lambda\\
&\leq&
C~\sum_{2k+\ell= J+1}\omega_{k,\ell}~R^{Q(p-1)+p(2k+\frac{\ell}{2})}
\int_{0<|\lambda|\leq\gamma}\big((2|\alpha|+n)|\lambda|\big)^{p(
k+\frac{\ell}{2})-{\sigma}}~|\lambda|^{n}d\lambda\\
&\leq&
2C~\sum_{\ell=0}^{J+1}~\omega_{\ell}~R^{Q(p-1)+p(J+1-\frac{\ell}{2})}
\int_{0}^{\gamma}\big((2|\alpha|+n)|\lambda|\big)^
{p(\frac{J+1}{2})-{\sigma}}~
|\lambda|^{n}d\lambda.
\end{eqnarray*}}
That is
\begin{equation}\label{3.3.16}
 S_{1}\leq C~R^{Q(p-1)+p(J+1-\frac{\ell}{2})}
\gamma^{p(\frac{J+1}{2})+{\frac{Q}{2}}-\sigma},~~\forall
\ell=0,1,...,J+1,
\end{equation}
provided that $p(\frac{J+1}{2}) + \frac{Q}{2}-\sigma > 0,$ which follows from the inequality (\ref{3}).\\
Now to estimate $S_{2},$ we may apply H\"{o}lder's inequality for
$q=\frac{2}{p}$ and Plancherel formula. Thus, we immediately obtain
\begin{eqnarray*}
S_{2} &\leq&
\Bigg(\int_{\mathbb{R}}(\|\mathcal{F}(a_{k})(\lambda)\|^{p})^{\frac{2}{p}}
|\lambda|^{n}d\lambda\Bigg)^\frac{p}{2}
\Bigg(\int_{|\lambda|>\gamma}\big((2|\alpha|+n)|\lambda|\big)
^{\frac{2\sigma}{p-2}}|\lambda|^{n}d\lambda\Bigg)^{\frac{2-p}{2}}\\
&\leq & C \|\mathcal{F}(a_{k})\|^{p}_{\mathcal{L}^{2}}
\Bigg(\int_{|\lambda|>\gamma}\big((2|\alpha|+n)|\lambda|\big)
^{\frac{2\sigma}{p-2}}|\lambda|^{n}d\lambda\Bigg)^{\frac{2-p}{2}}\\
&\leq & 2 C \|\mathcal{F}(a_{k})\|^{p}_{\mathcal{L}^{2}}
\Bigg(\int_{\gamma}^{\infty}
\big((2|\alpha|+n)|\lambda|\big)^
{\frac{2\sigma}{p-2}}|\lambda|^{n}d\lambda\Bigg)^{\frac{2-p}{2}}\\
&\leq &C \|\mathcal{F}(a_{k})\|^{p}_{\mathcal{L}^{2}}
\gamma^{{Q\over4}(2-p)-\sigma}
\end{eqnarray*}
provided that ${Q\over4}(2-p)-\sigma<0,$ which is a consequence of the left hand side of (\ref{3}).
Thanks to Plancherel's formula for Laguerre Fourier  transform it
follows
\begin{eqnarray*}
 \|\mathcal{F}(a_{k})\|^{2}_{\mathcal{L}^{2}}
=\|a_{k}\|_{L^{2}(\mathbb{H}^{n})}^{2}
&=&\int_{\mathbb{H}^{n}}|a_{k}(z,t)|^{2}~dV(z,t)\nonumber\\
&\leq&|B(0,R)|^{1-\frac{2}{p}}\nonumber\\
&\leq& C~ R^{-Q(\frac{2-p}{p})}.
\end{eqnarray*}
That is  $$ \|\mathcal{F}(a_{k})\|^{p}_{\mathcal{L}^{2}}\leq
C~ R^{-\frac{Q}{2}(2-p)},$$ and hence,
\begin{equation}\label{3.3.19}
S_{2}\leq C~
R^{-\frac{Q}{2}(2-p)}\gamma^{{{Q\over4}(2-p)}-\sigma}.
\end{equation}
However, to prove that $S_1 + S_2 \leq C,$ we shall discuss the cases $0 < R < 1$ and $R\geq 1.$
Hence, in order to deal with the case $0<R<1,$ we need more precise
estimates, so we consider the set $\Gamma_{\gamma} ;$ the collection of all numbers $\gamma$ satisfying
\begin{equation*}\label{3.3.17}
\Gamma_{\gamma} =\Big\{\gamma>0,~\frac{\frac{Q}{2}(2-p)}{{\frac{Q}{4}(2-p)}-\sigma}\log(R)\leq
\log(\gamma)\leq\frac{Q(1-p)-p(J+1)}{p(\frac{J+1}{2})+{\frac{Q}{2}}-\sigma}\log(R)\Big\}.
\end{equation*}
We mention that the collection $\Gamma_{\gamma}$ above is an nonempty set if and only if
\begin{equation*}\label{C}
\frac{\frac{Q}{2}(2-p)}{\frac{Q}{4}(2-p)-\sigma}\times
\frac{p(\frac{J+1}{2})+{\frac{Q}{2}}-\sigma}{Q(1-p)-p(J+1)}\leq1
\end{equation*}
which is a different formulation of the  hand side of
(\ref{3}), that is $\frac{Q}{2}(2-p)\leq\sigma.$

Now let us choose $\gamma \in\Gamma_{\gamma}$ and using the fact that
$\frac{Q}{2}+p\frac{(J+1)}{2}-\sigma>0$ together with
the right hand side of (\ref{3}) it follows that
\begin{equation}\label{3.3.18}
 S_{1}\leq C~R^{Q(p-1)+p(J+1)}
\gamma^{p(\frac{J+1}{2})+{\frac{Q}{2}}-\sigma}.
\end{equation}
Also, with the same choose of $\gamma \in\Gamma_{\gamma}$ and under the condition $\frac{Q}{2}(2 - p) < \sigma,$ together with
the help of the left hand side of (\ref{3}) we obtain
\begin{equation}\label{3.3.20}
S_{2}\leq C.
\end{equation}
Combining (\ref{3.3.18}) and (\ref{3.3.20}) we obtain
\begin{equation}\label{R<1}
S_{1}+S_{2}\leq C\qquad \mbox{for}\qquad0<R<1.
\end{equation}
Now, to deal with the case $ R\geq 1,$ we may take
\begin{equation}\label{3.3.21}
\gamma =R^{\frac{Q(1-p)-p(J+1)}{p(\frac{J+1}{2})+{\frac{Q}{2}}-\sigma}}
\end{equation}
so, using the fact that $ R \geq 1,$ we obtain
\begin{equation}\label{3.3.22}
\gamma \leq
 R^{\frac{\frac{Q}{2}(2-p)}{\frac{Q}{4}(2-p)-\sigma}}.
\end{equation}
which leads to
\begin{equation}\label{R>1}
S_{1}+S_{2}\leq C\qquad \mbox{for}\qquad R\geq1.
\end{equation}
Hence, to prove (\ref{3.3.15}), it is enough to combine (\ref{R<1}) and (\ref{R>1}). The proof of the main theorem is completed.

\textbf{Acknowledgements.}
\textit{This project was supported by King Saud University, Deanship of
Scientific Research, College of Science Research Center}.

%%%%%%%%%%%%%%%%%%%%%%%%%%%%%%%%%%%%%%%%%%%%%%%%%%%%%%%%%%%%%%%%%%%%%%%%%%%%
%%%%%%%%%%%%    bibliography     %%%%%%%%%%%%
%%%%%%%%%%%%%%%%%%%%%%%%%%%%%%%%%%%%%%%%%%%%%%%%%%%%%%%%%%%%%%%%%%%%%%%%%%%%


\begin{thebibliography}{20}


\bibitem{AM1} {M. Assal}, {\it Hardy's type inequality associated with the Hankel transform for overcritical exponent}, Integr. Transf. Spec. F.,
(2010), 1-6.
\bibitem{AR} { M. Assal and A. Rahmouni}, {\it Hardy's type inequality associated with the Laguerre Fourier transform}, Integr. Transf. Spec. F., 24, (2013), 156--163.

\bibitem{AR1} { M. Assal and A. Rahmouni}, {\it An improved Hardy's inequality associated with the Laguerre Fourier transform}, Collect. Math., (2013),  1--11.

\bibitem{CO} R. R. Coifman,
{\it A real-variable characterization of $H^p,$}  Studia Math. {\bf
51}, (1974), 269-274.


\bibitem{FS} C. Fefferman and E. M. Stein,
{\it $H^p$ spaces of several variables,} Acta Math. {\bf 129},
(1972), 137-193.


\bibitem{S3} G. B. Folland and E. M. Stein, \emph{Hardy Spaces on Homogeneous
Groups,} Princeton University Press, Princeton, NJ, 1982.

\bibitem{GR} J. Garcia-Cuerva and J. Rubio de Francia,
\emph{ Weighted Norm Inequalities and Related Topics,}  North
Holland, 1985.

\bibitem{Giulini} S. Giulini,
\emph{Bernstein and Jackson theorems for thé Heisenberg group,}
 J. Austral. Math. Soc. (Series A) {\bf 38} (1985), 241-254.

 \bibitem{Liu} H. P. Liu, {\it The group Fourier transforms and mulitipliers of the Hardy spaces on the Heisenberg group}, Approx. Theory $\&$ Its Appl.,
    {\bf 7} (1991), 106-117.
\bibitem{Lin} C. C. Lin, {\it $L^p$ multipliers and their $H^1 - L^1$ estimates on the Heisenberg group}, Revista Math. Ibero., {\bf 11} (1995), 269-308.

\bibitem{Lin1} C. C. Lin, {\it   H$\ddot{o}$rmander's $H^{p}$ Multiplier theorem for the Heisenberg group,}
     J. London Math. Soc. ({\bf3}) 67 (2003) 686-700.

\bibitem{kan2}{Y. Kanjin,} {\it On Hardy-Type Inequalities and Hankel Transforms,}
Monatshefte $f\ddot{u}$r Mathematik,  {\bf 127}, (1999), 311-319.

\bibitem{kan}{Y. Kanjin,} {\it Hardy's inequalities for Hermite and Laguerre expansions,} Bull. London Math. Soc.,  {\bf 29}, (1997), 331-337.


%\bibitem{RADHA} R. Radha, {\it Hardy type inequalities,}Taiwanese J. Math., {\bf 4}, (2000), 447-456.

\bibitem{RT} R. Radha and S. Thangavelu, {\it Hardy's inequalities for Hermite and Laguerre expansions,}
Proc. Amer. Math. Soc., {\bf 132}, (12), (2004), 3525-3536.

\bibitem{AR3} { A. Rahmouni and M. Assal},
{\it Hardy's type inequality for the  critical  exponent
associated with the inverse Laguerre Fourier transform,}
Integr. Transf. Spec. F., (2013) 1--9.
%\bibitem{Sababheh1} M. Sababheh, {\it Two-sided probabilistic versions of Hardy's inequality},
%J. Fourier An. Appl., {\bf 13}, (2007), 577-587.

%\bibitem{Sard}
%\textrm{A. Sard},
%\newblock{\em Linear Approximation,}
%\newblock Amer. Math. Soc., Providence, R.I., 1963.

%\bibitem{sat} M. Satake, {\it Hardy's inequalities for Laguerre expansions,} J. Math. Soc., {\bf 52}, (1), (2000), 17-24.

%\bibitem{Sababheh2} M. Sababheh,, {\it On an argument of K\"{o}rner and Hardy's inequality}, Analysis Mathematica, 34, (2008), 51-57.
%\bibitem{Sababheh3} M. Sababheh,  {\it A Study Of The Real Hardy Inequality} Journal of inequalities in pure and applied mathematics,
%{\bf 10}, Issue 4, (2009), [ONLINE
%http://jipam.vu.edu.au/article.php?sid=1160].
% \bibitem{Sababheh4} M. Sababheh,, {\it Hardy Inequalities On The Real Line}, Canadian Mathematical Bulletin,  to appear, (2010).


%\bibitem{Stroud} \textrm{A. H. Stroud},
%\newblock{\em Approximate calculation of multiple integrals,}
%\newblock Englewood Cliffs, N.J. Prentice-Hall, 1971.


%\bibitem{S3} E. M. Stein, \emph{Harmonic Analysis,} Princeton Univ. Press, Princeton, NJ, 1993.

\bibitem{S2} E. M. Stein, \emph{Harmonic Analysis, real variable Methods, orthogonality and
oscillatory integrals,} Princeton Univ. Press, Princeton, NJ, 1993.


\bibitem{TW} M. H. Taibleson and G. Weiss,
{\it The molecular characterization of certain Hardy spaces,}
Ast\'erisque {\bf 77}, (1980), Soci\'et\'e Math. de France, Paris,
67-149.

\bibitem{thang} S. Thangavelu, {\it On regularity of twisted spherical means and special Hermite
expansion}, Proc. Ind. Acad. Sci., {\bf 103}, (1993), 303-320.



\end{thebibliography}
\end{document}